\documentclass{amsart}
\usepackage{amssymb,latexsym}
\numberwithin{equation}{section}
\newtheorem{theorem}{Theorem}[section]

\newtheorem{corollary}[theorem]{Corollary}

\newtheorem{remark}[theorem]{Remark}

\newtheorem{example}[theorem]{Example}

\begin{document}

\pagenumbering{arabic} \pagestyle{headings}
\def\sof{\hfill\rule{2mm}{2mm}}
\def\ls{\leqslant}
\def\gs{\geqslant}
\def\SS{\frak S}
\def\T{\mathcal{T}}
\def\qq{{\bold q}}
\def\txx{{\tfrac1{2\sqrt{x}}}}

\title{ Permutations containing a pattern exactly
        once and avoiding at least two patterns of three letters}

\author{T. Mansour}
\maketitle
\begin{center}
LaBRI (UMR 5800), Universit\'e Bordeaux 1, 351 cours de la Lib\'eration\\
       33405 Talence Cedex, France\\[4pt]
{\tt toufik@labri.fr}
\end{center}
\markboth{Toufik Mansour}{Restricted certain patterns}

\section*{Abstract}
In this paper, we find an explicit formulas, or recurrences, in
terms of generating functions for the cardinalities of the sets
$S_n(T;\tau)$ of all permutations in $S_n$ that contain $\tau\in
S_k$ exactly once and avoid a subset $T\subseteq S_3$, $|T|\geq2$.
The main body of the paper is divided into three sections
corresponding to the cases $|T|=2,3$ and $|T|\geq4$.
%====================================================================
\section{Introduction}
Let $[p]=\{1,\dots,p\}$ denote a totally ordered alphabet on $p$
letters, and let $\alpha=(\alpha_1,\dots,\alpha_m)\in [p_1]^m$,
$\beta=(\beta_1,\dots,\beta_m)\in [p_2]^m$. We say that $\alpha$
is {\it order-isomorphic\/} to $\beta$ if for all $1\leq i<j\leq
m$ one has $\alpha_i<\alpha_j$ if and only if $\beta_i<\beta_j$.
For two permutations $\pi\in S_n$ and $\tau\in S_k$, an {\it
occurrence\/} of $\tau$ in $\pi$ is a subsequence $1\leq
i_1<i_2<\dots<i_k\leq n$ such that $(\pi_{i_1}, \dots,\pi_{i_k})$
is order-isomorphic to $\tau$; in such a context $\tau$ is usually
called the {\it pattern\/}. We say that $\pi$ {\it avoids\/}
(respectively, {\em contains exactly once}) $\tau$, or is
$\tau$-{\it avoiding\/}, if there is no occurrence of $\tau$ in
$\pi$ (respectively, there is exactly one occurrence of $\tau$ in
$\pi$). Pattern avoidance proved to be a useful language in a
variety of seemingly unrelated problems, from stack sorting
\cite[Ch. 2.2.1]{Kn} to singularities of Schubert varieties
\cite{LS}. A natural generalization of single pattern avoidance is
{\it subset avoidance\/}; that is, we say that $\pi\in S_n$ avoids
a subset $T\subset S_k$ if
$\pi$ avoids any $\tau\in T$.\\

Two sets, $T_1$, $T_2$, are said to be {\em Wilf equivalent} (or to belong to
the same {\em Wilf class}) if and only if $|S_n(T_1)|=|S_n(T_2)|$ for
any $n\geq 0$; the Wilf class of $T$ we denote by $\overline{T}$. Also
for two sets $(T_1;\tau^1)$, $(T_2;\tau^2)$, are said to belong
to the same {\em almost Wilf class} if and only if $|S_n(T_1;\tau^1)|=|S_n(T_2;\tau^2)|$
for any $n\geq 0$.\\

The study of the sets $S_n(\tau)$ was initiated by Knuth
~\cite{Kn}, who proved that $|S_n(\tau)|=\frac{1}{n+1}{{2n}\choose
n}$ for any $\tau\in S_3$. Knuth's results where further extended
in two directions. West ~\cite{W} and Stankova ~\cite{St} analyzed
$S_n(\tau)$ for $\tau\in S_4$ and obtained the complete
classification, which contains $3$ distinct Wilf classes. On the
other hand, Simion and Schmidt ~\cite{SS} studied $S_n(T)$ for
arbitrary subsets $T\subseteq S_3$ and discovered $7$ Wilf
classes. The study of $S_n(\alpha,\tau)$ for all $\alpha\in S_3$,
$\tau\in S_4(\alpha)$, was completed by West ~\cite{W}; Billey,
Jockusch and Stanley ~\cite{B};
and Guibert ~\cite{G}.\\

Later than, Mansour and Vainshtein \cite{MV2,MV3,MV4} presented a
recurrence formula to the generating function for the number of
$132$-avoiding permutations and avoiding or containing exactly
once another pattern in $S_k$ (see \cite{CW,Kr}). Besides, Mansour
\cite{M} gave a complete answer for avoiding at least two pattern
in $S_3$ and avoiding another pattern in $S_k$. On the other hand,
Robertson ~\cite{R} presented the complete classification of the
almost Wilf classes for avoiding or containing exactly once a
pattern from $S_3$, and containing exactly once another pattern
from $S_3$.\\

In the present paper, we find an explicit formulas, or
recurrences, in terms of generating functions for the
cardinalities of the sets $S_n(T;\tau)$ of all permutations in
$S_n$ that contain $\tau\in S_k$ exactly once and avoid a set $T$,
$|T|\geq2$, of patterns from $S_3$. Throughout the paper, we often
make use of the following remark.

\begin{remark}
\label{rem1}
In \cite{W} observed that if $\pi\in S_k$ contain any pattern $\tau\in T$,
then $S_n(T,\pi)=S_n(T)$. Besides, $|S_n(T;\pi)|=0$ for all $n\geq k$.
On the other hand, by Erd{\"o}s and Szekeres \cite{ES} we obtain that
$|S_n(T;\tau)|=0$ for all $n\geq 7$ where $\{123,321\}\subset T$.
Since that, from now we suppose that $\pi$ is $T$-avoiding,
and $\{123,321\}\not\subseteq T$.
\end{remark}

The main body of the paper is divided into three sections
corresponding to the cases $|T|=2,3$ and $|T|\geq4$.
%======================================================================
\section{A pair}
In this section we present, by explicit formula or recurrence in
terms of generating functions, the cardinality of the sets
$S_n(\beta,\gamma;\tau)$ where $\beta,\gamma\in S_3$, $\tau\in
S_k$, $k\geq 3$. By the three natural operations, the complement,
the reversal, the inverse (see Simion and Schmidt
\cite[Lem.~1]{SS}), and Remark \ref{rem1} we have to consider the
following four possibilities:

$$\begin{array}{rlll}
  1)&   S_n(132,123;\tau),  & \ \mbox{where}    & \ \tau\in S_k(132,123), \\
  2)&   S_n(132,321;\tau),  & \ \mbox{where}    & \ \tau\in S_k(132,321), \\
  3)&   S_n(132,213;\tau),  & \ \mbox{where}    & \ \tau\in S_k(132,213), \\
  4)&   S_n(132,231;\tau),  & \ \mbox{where}    & \ \tau\in S_k(132,231).
\end{array}$$

The main body of this section is divided into four subsections corresponding
to the above four cases.
%==========2=1===putk5=========
\subsection{$\bf T=\{123, 132\}$.}

Let $a_\tau(n)=|S_n(123,132;\tau)|$, and
$a'_\tau(n)=|S_n(123,132,\tau)|$; and let $A_\tau(x)$
($A'_\tau(x)$) be the ordinary generating function for the
sequence $\{a_\tau(n)\}_{n\geq0}$ ($\{a'_\tau(n)\}_{n\geq 0}$). We
find a recurrence for the generating function $A_\tau(x)$.

\begin{theorem}
\label{a12}
Let $k\geq 2$, and $\tau\in S_k(132,123)$; so there exist $r$ maximal such that
$\tau=(k-1,k-2,\dots,k-r+1,k,\tau')$. Then

{\rm (i)}  for $r\geq 2$,
        $$A_{(k-1,k-2,\dots,k-r+1,k,\tau')}(x)=\left\{
        \begin{array}{ll}
        \frac{x^r(1-x)}{1-2x+x^r} A_{\tau'}(x)      &\ ;\tau'\neq\emptyset \\
        \frac{x^r(1-x)}{1-2x+x^r} A'_{\tau}(x)      &\ ;\tau'=\emptyset.
        \end{array} \right. .$$

{\rm (ii)} for $m$ maximal such that
$\tau=(k,k-1,\dots,k-m+1,\tau')$, {\small
$$A_{\tau}(x)=\left\{
        \begin{array}{ll}
        xA_{(k-1,\dots,k-m+1,\tau')}(x)+\sum\limits_{j=2}^{m} x^{j+1}A_{(k-j,\dots,k-m+1,\tau')}(x) &\ ;m\geq 2 \\
        (x-x^2-x^3-\dots-x^{k-1}) A_{\tau'}(x) &\ ;m=1
        \end{array} \right. .$$}
Besides, $A_{(1)}(x)=x$.
\end{theorem}
\begin{proof}
Let $\tau\in S_k(132,123)$; by \cite[Th.~3(i)]{M} there exist $r$
maximal such that $\tau=(k-1,k-2,\dots,k-r+1,k,\tau')$. Let
$\alpha\in S_n(132,123)$; choose $t$ such that $\alpha_t=n$, so by
\cite[Th.~3(i)]{M} $\alpha=(n-1,\dots,n-t+1,n,\alpha')$.
\begin{enumerate}
\item   Let $r\geq 2$ and $\tau'\neq\emptyset$. By definitions it is easy to obtain
    for any $n\geq k$,
        $$a_\tau(n)=\sum_{t=1}^{r-1} a_\tau(n-t) +a_{\tau'}(n-r).$$

\item   Let $r\geq 2$ and $\tau'=\emptyset$. By definitions it is easy to get for $n\geq k$
        $$a_\tau(n)=\sum_{t=1}^{r-1} a_{\tau}(n-t) +a'_\tau(n-r).$$

\item   For $r=2$, by use \cite[Th.~3(i)]{M} there exist $m$ maximal such that
    $\tau=(k,k-1,\dots,k-m+1,\tau')$. This yields
        $$a_\tau(n)=a_{(k-1,\dots,k-m+1,\tau')}(n-1)+ \sum_{j=2}^m a_{(k-j,\dots,k-m+1,\tau')}(n-1-j),$$
    for $m\geq 2$, and
        $$a_\tau(n)=a_{\tau'}(n-1)-\sum_{j=2}^{k-1} a_{\tau'}(n-j),$$
    where $m=1$.
\end{enumerate}
Besides, in the above cases $a_\tau(n)=0$ for all $n\leq k-1$ and
$a_\tau(k)=1$, hence in terms of generating function the theorem
holds.
\end{proof}

\begin{corollary}
Let $r_i\geq 2$ such that $r_1+\dots+r_m=k$, and
$\tau=(p_1,p_2,\dots,p_m)\in S_k(123,132)$ where
$p_i=(t_i-1,t_i-2,\dots,t_i-r_i+1)$, and $t_i=k-(r_1+\dots+r_{i-1})$
for $i=1,2,\dots,m$. Then
    $$A_{\tau}(x)=\frac{1-x}{1-2x+x^{r_m}}
        \prod_{i=1}^m \frac{x^{r_i}(1-x)}{1-2x+x^{r_i}}.$$
\end{corollary}
\begin{proof}
By induction and part $(i)$ of Theorem \ref{a12} we obtain that
    $$A_{\tau}(x)=A'_{p_m}(x)
        \prod_{i=1}^m \frac{x^{r_i}(1-x)}{1-2x+x^{r_i}},$$
so by \cite[Th.~3(ii)]{M} the corollary holds.
\end{proof}

\begin{example}\label{ex1}
For $k=3$, Theorem \ref{a12} and \cite[Th. 3(ii)]{M} yields
    $$A_{213}(x)=\frac{x^3}{(1-x-x^2)^2},\quad
      A_{231}(x)=A_{312}(x)=\frac{x^3}{1-x},\quad
      A_{321}(x)=x^3+3x^4.$$
\end{example}
%=========2=2===putk5=========
\subsection{$\bf T=\{132, 321\}$.}
Let $b_\tau(n)=|S_n(132,321;\tau)|$ and let $B_\tau(x)$ be the ordinary
generating functions of the sequences $\{b_\tau(n)\}_{n\geq0}$.
We find a explicit formula for the generating function $B_\tau(x)$.

\begin{theorem}\label{a26}
Let $k\geq 1$; then

{\rm (i)}
$$B_{(1,2,\dots,k)}(x)=x^k +2\sum\limits_{j=k+1}^{2k-1}
(2k-j)x^j;$$

{\rm (ii)} for $1\leq d\leq k-1$,
        $$B_{(d+1,d+2,\dots,k,1,2,\dots,d)}(x)=\frac{x^k}{1-x};$$

{\rm (iii)}    for $1\leq d\leq m-2\leq k-2$,
        $$B_{(d+1,d+2,\dots,m-1,1,2,\dots,d,m,m+1,\dots,k)}(x)=x^k.$$
\end{theorem}
\begin{proof}
Let $\tau\in S_k(132, 321)$; by \cite[Th.~6(i)]{M} there exist
$m$, $2\leq m\leq k+1$, and $d$, $1\leq d\leq m-2$ such that
    $$\tau=(d+1,d+2,\dots,m-1,1,2,\dots,d,m,m+1,\dots,k).$$
Similarly, for $\alpha\in S_n(132,321)$ there exist $r$, $2\leq r\leq n+1$,
and $t$, $1\leq t\leq r-2$ such that
    $$\alpha=(t+1,t+2,\dots,r-1,1,2,\dots,t,r,r+1,\dots,n).$$
By this fact, the rest it is easy to see. \end{proof}

\begin{example}
\label{ex2}
Theorem \ref{a26} yields,
    $$B_{123}(x)=x^3+4x^4+2x^5,\quad
      B_{213}(x)=x^3,\quad
      B_{231}(x)=B_{312}(x)=\frac{x^3}{1-x}.$$
\end{example}
%==========2=3===putk5===========
\subsection{$\bf T=\{132, 213\}$.}
Let $c_\tau(n)=|S_n(132,213;\tau)|$ and let $C_\tau(x)$ be the ordinary
generating functions of the sequences $\{c_\tau(n)\}_{n\geq0}$.
We find a recurrence formula for the generating function $C_\tau(x)$.

\begin{theorem}\label{a23}
Let $\tau\in S_k(132,213)$; then

{\rm(i)}  there exist $k+1=r_0>r_1>\dots>r_m=1$ such that
        $$\tau=(r_1,r_1+1,\dots,k,r_2,r_2+1,\dots,r_1-1,\dots,r_m,r_m+1,\dots,r_{m-1}-1);$$

{\rm (ii)} for all $0\leq r\leq k-1$,
        $$C_{(r+1,\dots,k,\tau)}(x)=\frac{x^{k-r}(1-x)}{1-2x+x^{k-r}} C_{\tau'}(x),$$
        where $\tau'\neq\emptyset$, and
        $$C_{(1,2,\dots,k)}(x)=\frac{x^k(1-x)^2}{(1-2x+x^k)^2}.$$
\end{theorem}
\begin{proof}
$(i)$ holds by \cite[Th.~8(i)]{M}. To prove $(ii)$, let $\alpha\in
S_n(132,213;\tau)$, so by $(i)$ there exist
$n+1=t_0>t_1>\dots>t_m\geq 1$ such that
    $$\alpha=(t_1,t_1+1,\dots,t_0-1,t_2,t_2+1,\dots,t_1-1,\dots,t_m,t_m+1,\dots,t_{m-1}-1),$$
therefore for $\tau'\neq\emptyset$
    $$c_\tau(n)=\sum_{j=n-k+r_1+1}^n c_\tau(j-1)+c_{\tau'}(n-k+r_1-1).$$
When $\tau'=\emptyset$ which means that $\tau=(1,2,\dots,k)$ we get that
    $$c_{\tau}(n)=\sum_{j=n-k+2}^n c_{\tau}(j-1)+c'_{\tau}(n),$$
where $c'_{\tau}(n)$ is the number of all permutations
$\{132,213\}$-avoiding and $\tau$-avoiding. Besides, $c_\tau(k)=1$
and $c_\tau(n)=0$ for $n\leq k-1$, hence the rest is easy to see
by use \cite[Th.~8(ii)]{M}.
\end{proof}

\begin{corollary}
Let $k\geq 1$; then $C_{(k,k-1,\dots,1)}(x)=x^k$.
\end{corollary}

\begin{example}\label{ex3}
Theorem \ref{a23} yields
$$C_{123}(x)=\frac{x^3}{(1-x-x^2)^2},\quad C_{231}(x)=C_{312}(x)=\frac{x^3}{1-x},\quad C_{321}(x)=x^3.$$
\end{example}
%============2=4===putk5=============
\subsection{$\bf T=\{213, 231\}$.}
Let $d_\tau(n)=|S_n(213,231;\tau)|$ and let $D_\tau(x)$ be the
ordinary generating functions of the sequences
$\{d_\tau(n)\}_{n\geq0}$. We find a recurrence formula for the
generating function $D_\tau(x)$.\\

Let $\tau\in S_k(132,231)$; by \cite[Th.~11]{M} $\tau_i$ is either
the right maximum, or the right minimum, for all $1\leq i\leq
k-1$.

\begin{theorem}
\label{a24} Let $\tau\in S_k(132,231)$; then

{\rm (i)}  there exist $\tau'$ such that either
        $\tau=(k,k-1,\dots,k-r+1,\tau',k-r)$, or
        $\tau=(k-r,\tau',k-r+1,\dots,k)$, where $1\leq r\leq k-1$;

{\rm (ii)} $$D_{\tau}(x)=\frac{x^{r+1}}{(1-x)^r} D_{\tau'}(x),$$
        where $\tau'\neq\emptyset$, and
        $$D_{(1,2,\dots,k)}(x)=D_{(k,\dots,2,1)}(x)=\frac{x^k}{(1-x)^{k-1}}.$$
\end{theorem}
\begin{proof}
If $\tau_1\neq k$ and $\tau_k\neq k$, then $\tau$ contain $132$ or $231$, so
$\tau_1=k$ or $\tau_k=k$. By use induction on $k$ the first case holds.

Let $\alpha\in S_n(132,231)$; similarly to $(i)$ there exist $\alpha'$ such that
$\alpha=(n,\dots,n-t+1,\alpha',n-t)$, or $\alpha=(n-t,\alpha',n-t+1,\dots,n)$,
where $1\leq t\leq n-1$. Then for $0\leq m\leq r-1$,
    $$d_{(k-m,\dots,k-r+1,\tau',k-r)}(n)=d_{(k-m,\dots,k-r+1,\tau',k-r)}(n-1)+d_{(k-1-m,\dots,k-r+1,\tau',k-r)}(n-1),$$
and
    $$d_{(k-r,\tau',k-r+1,\dots,k-m)}(n)=d_{(k-r,\tau',k-r+1,\dots,k-m)}(n-1)+d_{(k-r,\tau',k-r+1,\dots,k-1-m)}(n-1).$$
Besides $d_\tau(n)=0$ for all $n\leq k-1$ and $d_\tau(k)=1$, hence it is easy to obtain
    $$D_\tau(x)=\frac{x^{r+1}}{(1-x)^r} D_{\tau'}(x).$$
The rest is obtain immediately by above recurrence.
\end{proof}

Let us denote the sequence $k\dots (r+2)(r+1)$ by $<k,r>$.

\begin{corollary}
Let $k\geq 1$, and $\tau=(<k,r_1>,<r_1-1,r_2>,\dots,<r_{m-1}-1,r_m>,r_{m-1},r_{m-2},\dots,r_1)\in S_k(132,231)$.
Then
    $$D_{\tau}(x)=\frac{x^k}{(1-x)^{k-m}}.$$
\end{corollary}
\begin{proof}
By Theorem \ref{a24} we obtain that
    $$D_{\tau}(x)=\prod_{i=1}^{m-1} \frac{x^{r_{i-1}-r_i}+1}{(1-x)^{r_{i-1}-r_i}} D_{(r_m,\dots ,2,1)}(x),$$
which means that
    $$D_{\tau}(x)=\prod_{i=1}^{m-1} \frac{x^{r_{i-1}-r_i}+1}{(1-x)^{r_{i-1}-r_i}} \cdot \frac{x^{r_m}}{(1-x)^{r_m-1}},$$
hence the corollary holds.
\end{proof}

\begin{example}
\label{ex4}
Theorem \ref{a24} yields,
$$D_{123}(x)=D_{321}(x)=\frac{x^3}{(1-x)^2},\quad D_{213}(x)=D_{312}(x)=\frac{x^3}{1-x}.$$
\end{example}
%====================================================================
%====================================================================
\section{A triplet}

In this section, we calculate the cardinality of the sets
$S_n(T;\tau)$ such that $T\subset S_3$, $|T|=3$ and $\tau\in
S_k(T)$ for $k\geq 3$. By three natural operations the complement,
the reversal, the inverse (see Simion and Schmidt
\cite[Lem.~1]{SS}), and Remark \ref{rem1} we have to consider the
following four possibilities:

$$\begin{array}{rlll}
  1)&   S_n(123,132,213;\tau),  & \ \mbox{where}    & \ \tau\in S_k(123,132,213), \\
  2)&   S_n(123,132,231;\tau),  & \ \mbox{where}    & \ \tau\in S_k(123,132,231), \\
  3)&   S_n(123,231,312;\tau),  & \ \mbox{where}    & \ \tau\in S_k(123,231,312), \\
  4)&   S_n(132,213,231;\tau),  & \ \mbox{where}    & \ \tau\in S_k(132,213,231).
\end{array}$$

The main body of this section is divided into four subsections corresponding
to the above four cases.
%===================1===putk345============
\subsection{$\bf T=\{123, 132, 213\}$.}
Let $e_\tau(n)=|S_n(123,132,213;\tau)|$ and let $E_\tau(x)$ be the ordinary
generating functions of the sequences $\{e_\tau(n)\}_{n\geq0}$.
We find a recurrence formula for the generating function $E_\tau(x)$.

\begin{theorem}
\label{a123}
Let $k\geq 4$, $\tau\in S_k(123,132,213)$. Then
$$E_{(k-1,k,\tau')}(x)=\frac{x^2}{1-x} E_{\tau'}(x),\quad E_{(k,\tau')}(x)=xE_{\tau'}(x),$$
Besides $E_{\tau}(x)$ is given by $x^4$, $\frac{x^3}{1-x}$,
$\frac{x^3}{1-x}$, $x^3$, $x^2$, $\frac{x^2}{(1-x)^2}$, $x$, where
$\tau=4231$, $231$, $312$, $321$, $21$, $12$, $1$; respectively.
\end{theorem}
\begin{proof}
Let $\alpha\in S_n(123,132,213;\tau)$; by \cite[Lem.~14]{M} we
obtain immediately $\alpha=(n-1,n,\alpha')$ or
$\alpha=(n,\alpha')$.

Now let $k\geq 4$, and $\tau=(k-1,k,\tau')$; so for all $n\geq k$
    $$e_\tau(n)=d_{\tau}(n-1)+d_{\tau'}(n-2),$$
which means the first case holds.

Now let $k\geq 4$, and $\tau=(k,\tau')\neq 4231$; so for all $n\geq 4$
    $$e_\tau(n)=e_{\tau'}(n-1),$$
which means the second case holds. Similarly we have all the
special cases.
\end{proof}

\begin{example}
By Theorem \ref{a123} we have immediately for all $k\geq 4$ that
    $$E_{(k,\dots,2,1)}(x)=E_{(k,\dots,4,2,3,1)}(x)=x^k.$$
Another example, for $k\geq 2$
    $$E_{(k-1,k,k-3,k-2,\dots,1,2)}(x)=
        \left\{ \begin{array}{ll} \frac{x^k}{(1-x)^{k/2+1}}, & \ \ k=2,4,6,\dots\\
                     \frac{x^k}{(1-x)^{(k-1)/2}}, &\  \ k=3,5,7,\dots
        \end{array}\right..$$
\end{example}
%=================2===putk5===============
\subsection{$\bf T=\{123, 132, 231\}$.}
Let $f_\tau(n)=|S_n(123,132,231;\tau)|$ and let $F_\tau(x)$ be the
ordinary generating functions of the sequences
$\{f_\tau(n)\}_{n\geq0}$. We find an explicit formula for the
generating function $F_\tau(x)$.

Let $\tau\in S_k(123,132,231)$; by \cite[Th.~17]{M} we have
    $$\tau=(k,k-1,\dots,k-r+1,k-r-1,\dots,1,k-r),$$
where $1\leq r\leq k$. Therefore

\begin{theorem}
Let $k\geq 3$, and $k-2\geq r\geq 1$; then
\label{a124}
$$\begin{array}{l}
F_{(k,\dots,21,)}(x)=x^k+(k-1)x^{k+1};\\
F_{(k-1,\dots,2,1,k)}(x)=\frac{x^k}{1-x};\\
F_{(k,\dots,k-r+1,k-r-1,\dots,1,k-r)}(x)=x^k.
\end{array}$$
\end{theorem}
%=================3===putk5===============
\subsection{$\bf T=\{123, 231, 312\}$.}
Let $g_\tau(n)=|S_n(123,231,312;\tau)|$ and let $G_\tau(x)$ be the
ordinary generating functions of the sequences
$\{g_\tau(n)\}_{n\geq0}$. We find a recurrence formula for the
generating function $G_\tau(x)$.

Let $\tau\in S_k(123,231,312)$; by \cite[Th.~21]{M} we get
    $$\tau=(r,r-1,\dots,1,k,k-1,\dots,r+1),$$
where $1\leq r\leq k$. Therefore

\begin{theorem}
Let $k\geq 3$, and $k-1\geq r\geq 1$; then
\label{a145}
$$\begin{array}{l}
G_{(k,\dots,2,1}(x)=\frac{x^k(1+x)}{1-x};\\
G_{(r,r-1,\dots,1,k,k-1,\dots,r+1)}(x)=x^k.
\end{array}$$
\end{theorem}
%=================4===putk5===============
\subsection{$\bf T=\{132, 213, 231\}$.}
Let $h_\tau(n)=|S_n(132,213,231;\tau)|$ and let $H_\tau(x)$ be the
ordinary generating functions of the sequences
$\{h_\tau(n)\}_{n\geq0}$. We find a recurrence formula for the
generating function $H_\tau(x)$.

Let $\tau\in S_k(132,213,231)$; by \cite[Th.~23]{M} we get
        $$\tau=(k,k-1,\dots,r+1,1,2,\dots,r),$$
where $1\leq r\leq k$. Therefore

\begin{theorem}
Let $k\geq 3$, and $k-1\geq r\geq 1$; then
\label{a234}
$$\begin{array}{l}
H_{(1,2,\dots,k)}(x)=\frac{x^k}{1-x};\\
H_{(k,k-1,\dots,r+1,1,2,\dots,r)}(x)=x^k.
\end{array}$$
\end{theorem}
%=========================================================================
\section{A quartet and A quintet}
By Erd{\"o}s ~\cite{ES} (see also Simion and Schmidt
\cite[pr.~17]{SS}) $|S_n(T)|=0$ for all subsets $T$ of $S_3$ such
that $\{123,321\}\subseteq T$, and on the other hand
$|S_n(T)|=2,1$ for all $\{123,321\}\not\subset T\subseteq S_3$
such that $|T|=4,5$. All these yields the following theorem.

\begin{theorem}
\label{a1234}
Let $\tau\in S_k$. Then
$$\begin{array}{ll}
|S_n(123,132,213,231;\tau)|&=\left\{ \begin{array}{ll}
                    1,&\ \mbox{where}\ n=k,\ \tau=(k,\dots,3,2,1), (k,\dots,3,1,2)\\
                    0,&\ \mbox{otherwise}
                    \end{array} \right.; \\
               & \\
|S_n(123,132,231,312;\tau)|&=\left\{ \begin{array}{ll}
                    1,&\ \mbox{where}\ n=k,k+1,\ \tau=(k,\dots,2,1)\\
                    1,&\ \mbox{where}\ n=k,\ \tau=(k-1,\dots,2,1,k)\\
                    0,&\ \mbox{otherwise}
                    \end{array} \right.;\\
                       & \\
|S_n(132,213,231,312;\tau)|&=\left\{ \begin{array}{ll}
                    1,&\ \mbox{where}\ n=k,\ \tau=(k,\dots,2,1), (1,2,\dots,k)\\
                    0,&\ \mbox{otherwise}
                    \end{array} \right.; \\
                &\\
|S_n(S_3\backslash\{123\};\tau)|&=\delta_{n,k}\delta_{\tau,(1,2,\dots,k)};\\
               &\\
|S_n(S_3;\tau)|&=0.
\end{array}$$
\end{theorem}
%==========================================================================

\end{document}